\documentstyle[11pt]{article}

\begin{document}

\input prepictex
\input pictex
\input postpictex

\def\R{{\rm I\!R}} 
\def\P{{\rm I\!P}} 
\def\N{{\rm I\!N}} 
\def\F{{\rm I\!F}}
\def\Q{{\mathchoice {\setbox0=\hbox{$\displaystyle\rm Q$}\hbox{\raise
0.15\ht0\hbox to0pt{\kern0.4\wd0\vrule height0.8\ht0\hss}\box0}}
{\setbox0=\hbox{$\textstyle\rm Q$}\hbox{\raise
0.15\ht0\hbox to0pt{\kern0.4\wd0\vrule height0.8\ht0\hss}\box0}}
{\setbox0=\hbox{$\scriptstyle\rm Q$}\hbox{\raise
0.15\ht0\hbox to0pt{\kern0.4\wd0\vrule height0.7\ht0\hss}\box0}}
{\setbox0=\hbox{$\scriptscriptstyle\rm Q$}\hbox{\raise
0.15\ht0\hbox to0pt{\kern0.4\wd0\vrule height0.7\ht0\hss}\box0}}}}
\def\Z{{\mathchoice {\hbox{$\sf\textstyle Z\kern-0.4em Z$}}
{\hbox{$\sf\textstyle Z\kern-0.4em Z$}}
{\hbox{$\sf\scriptstyle Z\kern-0.3em Z$}}
{\hbox{$\sf\scriptscriptstyle Z\kern-0.2em Z$}}}}

\def\OO{{\cal O}}
\def\LL{{\cal L}}
\def\MM{{\cal M}}

\newcommand{\vv}{\vspace{0.3cm}} 
\newcommand{\vV}{\vspace{0.5cm}} 
\newcommand{\VV}{\vspace{2cm}}

\title{Even Sets of Lines on  Quartic Surfaces}
\author{Wolf P. Barth}
\maketitle

\begin{center}
Mathematisches Institut der Universit\"{a}t

Erlangen, Bismarckstra\ss e $1\frac{1}{2}$

{e-mail: barth@mi.uni-erlangen.de}

\end{center}

\tableofcontents

\VV
\section{Introduction}

An effective divisor $D$ on a smooth (compact complex) surface $X$
is called {\em even}, if its class $[D] \in H^2(X,\Z)$
is divisible by $2$. $D$ may be assumed reduced w.l.o.g.
Then D being even is equivalent to the existence of 
a double cover $Y \to X$ branched exactly over $D$.

The aim of this note is to study arrangements $\Lambda(n)$ of 
$n \leq 10$ distinct lines
on a smooth quartic surface $X \subset \P_3$, which form an even
divisor in this sense. The result is that for $n \leq 8$ there are no
unexpected ones. To be precise, there are exactly the following five
even sets of $n \leq 8$ lines:

\begin{itemize}
\item[$\Lambda(6):$]  
{\em Two disjoint triangles.}
Here also three coplanar lines meeting in one point are
considered a triangle. 
\item[$\Lambda_1(8):$] {\em Eight mutually skew lines.} Such a union of
eight lines can be even, but need not always be it.
\item[$\Lambda_2(8):$] {\em A union of two  
disjoint space quadrangles.}
\item[$\Lambda_3(8):$] {\em The intersection of $X$ with a smooth 
quadric surface
consisting of four skew lines in each of the two rulings of the quadric.} 
\item[$\Lambda_4(8):$] {\em The intersection of $X$ with two planes,
in each plane a quadrangle.}
\end{itemize}

\noindent
The proof I give for this assertion is very crude: essentially checking cases.

And I give a partial classification for $n=10$ in the following sense:
Each even set $\Lambda(10)$ of ten lines on a smooth quartic surface
is of one of the types $\Lambda_1(10),...,\Lambda_{11}(10)$
given in section 4. Unfortunately at the moment I do not know
which of these types $\Lambda_i(10)$ do actually occur. However,
if a configuration $\Lambda_i(10)$ exists on a quartic surface,
it necessarily will be an even set of lines. The proof for this partial
classification is messy, again essentially checking cases. It doesn't
seem to make sense to pursue it further, say to even sets of twelve lines,
unless some new technique, adapted to this purpose evolves.

\VV
\section{Conditions}

Here I collect some conditions a divisor $\Lambda$ on a smooth
quartic $X \subset \R^3$ necessarily has to satisfy, if it consists of
lines.

\vv
$(\lambda 1):$ {\em The arrangement $\Lambda$ does not contain 
five lines in one plane.}

Of course, if five lines in $\Lambda$ would lie in one plane,
the surface $X$ would contain the plane and be singular.

\vv
$(\lambda 2):$ {\em If three lines from $\Lambda$ meet in one
point (concurrent lines), they will lie in a plane
(coplanar lines).}

Indeed, the plane is the tangent plane to $X$ at the point,
where the three lines meet.

This property $(\lambda 2)$ has the consequence, that 
we need
not distinguish between three lines in a plane forming a triangle
and three concurrent lines in a plane.

\vv
From now on I assume that the divisor
$$\Lambda = \sum L_i \subset X$$
is even, i.e., the class $[\Lambda] \in H^2(X,\Z)$ is divisible
by two. Say
$$\OO_X(\Lambda) = \LL^{\otimes 2}, \quad \LL \in Pic(X).$$ 
Then additionally we have the following properties.

\vv
$(\lambda 3):$ {\em The arrangement $\Lambda$ contains an even number $n$
of lines.}

Proof. If $\Lambda$ is even, its intersection number 
$(\Lambda.\P_2)=n$
with a general plane $\P_2$ must be even too.

\vv
$(\lambda 4):$ {Each line $L_i$ in $\Lambda$ meets an even number
$k_i$ of lines $L_j \not= L_i$.}

Proof. If $\Lambda$ is even, then each intersection number
$$(\Lambda.L_i)=-2+k_i$$
is even too.

\vv
$(\lambda 5):$ {\em Let $n$ be the number of lines in $\Lambda$. Then 
the integer 
$$k:=\frac{1}{2} \sum_i k_i$$
satisfies 
$$k \leq \frac{1}{2} \cdot n \cdot(n-2).$$ }
Indeed, $k_i \leq n-2$ for $i=1,...,n.$

\vv
$(\lambda 6):$ {\em There is the modulo-$4$ condition: 
$$k-n \mbox{ is divisible by }4.$$}
Proof. The self-intersection of the divisor $\Lambda$ is
$$\Lambda^2=-2 \cdot n + \sum_i k_i =2\cdot (k-n).$$
If now $\Lambda=2 \LL$, then 
$\Lambda^2=4 \cdot \LL^2 $
is divisible by $8$, because $\LL^2$ is even.

\vV
To recognize even sets of lines 
on a smooth quartic surface $X$ we use the following principles:

\vv
$(\pi 1):$ {\em If a set of disjoint smooth rational curves on a $K3$-surface
is even, then it contains either eight or 16 curves.}

This is lemma 3 of \cite{N}.

\vv
Let me call {\em elliptic fibre} any effective, reduced, connected
divisor $D \subset X$ with $D^2=0$ such that the linear system 
$|D|$ has no fixed component.  
Any elliptic fibre is a fibre
of some elliptic fibration $f:X \to \P_1$, so it is one of the 
reduced curves in Kodaira's table [BPV, p. 150] of elliptic fibres.

Proof. By Riemann-Roch $h^0(D) \geq 2$ and there is a fibration 
$f:X \to \P_1$ having $X$ as a fibre. The general fibre
of this fibration is smooth. Its connected components consist of
smooth elliptic curves $E$. The curve $E$ defines an elliptic 
pencil. The fibration $f$ is composed of this elliptic pencil.
Since $D$ is reduced and connected, it cannot be a multiple fibre, nor
consist of different fibres.
We have $[D]=[E]$ and $D$ indeed is a member of the elliptic pencil $|E|$.

\vv
$(\pi 2):$ {\em An elliptic fibre cannot be even.}

Proof. If $D$ were even, we would have $E=2F$ with $F$ another elliptic fibre.
The general member of $|F|$ would be smooth elliptic, and $E$ would consist of
two such components. As $E$ is reduced connected, this cannot happen.

\vv
$(\pi 3):$ {\em Let $D=D_1+D_2$ be the disjoint union
of two elliptic fibres $D_i$ on $X$. 
Then $[D_1]=[D_2]$ and in particular}
$$deg(D_1)=deg(D_2).$$ 

\vV

$(\pi 4):$ {\em Quadratic reduction:} Let $D=D_1+D_2$ be an arrangement of lines
on a quartic. Assume there is a quadric surface $Q \subset \P_3$
containing $D_1$, say 
$$Q.X = D_1+D_1'.$$
If $D$ is even, then also
$$D-Q.X = D_2-D_1'$$
is even, as well as $D_2+D_1'$. We say, the divisor $D_2+D_1'$ is
obtained from $D$ by quadratic reduction modulo $Q$.

\vv
$(\pi 5):$ {\em Elliptic reduction:} Let $D=D_1+D_2$ be an arrangement of lines
on a quartic with $D_1$ an elliptic fibre. Let $D_1'$ be any other
elliptic fibre in $|D_1|$. If $D$ is even, then also
$$D-(D_1+D_1')=D_2-D_1'$$
is even, as well as $D_2+D_1'$. We say, the divisor $D_2+D_1'$
is obtained from $D$ by elliptic reduction modulo $|D_1+D_1'|$.

\vV
We also recall the Riemann-Roch formula in two cases: 

For
the line bundle $\OO_X(m) \otimes \OO_X(-\Lambda)$:
\begin{eqnarray*}
\chi(\OO_X(m) \otimes \OO_X(-\Lambda))
 &=& 2 + \frac{1}{2}(4m^2-2mn+2(k-n)) \\
 &=& 2+2m^2+k-(m+1)n.
\end{eqnarray*}
This is $>0$ as soon as
$$k>(m+1)n-2(m^2+1).$$

And for the line bundle $\LL$ with $\LL^2=\frac{1}{2}(k-n)$ it says
$$\chi(\LL)=2+\frac{1}{2}\LL^2=2+\frac{1}{4}(k-n).$$
This is $>0$ as soon as $k>n-8$. Since $\LL^{\otimes 2} = \Lambda$
is effective, 
$$h^2(\LL)=h^0(-\LL)=0.$$
So $\LL$ itself is effective in this case.

\VV
\section{Examples ($n \leq 8$)}

First I want to show that the arrangements $\Lambda(6), \, 
\Lambda_1(8),...,\Lambda_4(8)$ indeed can be found on smooth quartic surfaces.

$\Lambda(6):$ Fix a line $L \subset \P_3$, intersection of two
planes $F_1$ and $F_2$. Fix triangles $L_1,L_2,L_3 \subset F_1$ and
$L_4,L_5,L_6 \subset F_2$ such that the six lines $L_1,...,L_6$
meet $L$ in six distinct points. Counting constants we find that
the vector space of quartic polynomials vanishing on the arrangement
$L_1,...,L_6$ has dimension $12$ at least. So the linear system
of quartics containing the arrangement has no base surface. Its general surface 
is irreducible.

Let $L_4',L_5',L_6'$ be three
lines in $F_1$ meeting $L$ precisely in its intersections
with $L_4,L_5,L_6$ respecively. Let $E_4,E_5,E_6$ be
the planes spanned by $L_4,L_4'$ etc. Then the six lines $L_1,...,L_6$
of our arrangement $\Lambda(6)$ are contained in the quartic 
$F_1+E_4+E_5+E_6$. This quartic is smooth on $F_1$ outside of 
$L_4',L_5',L_6'$. Bertini's theorem shows that the general
quartic in our linear system is smooth ootside of the 
three points $L \cap (L_4 \cup L_5 \cup L_6)$.  An analogous argument
shows that this quartic also is smooth away from
$L \cap(L_1 \cup L_2 \cup L_3)$. So the arrangement lies
on a smooth quartic surface.
  
\vv
$\Lambda_1(8)$: It is known that there are smooth quartics
$X \subset \P_3$ carrying $16$ disjoint lines. (Their
equations are computed in \cite{BN}, but they were known already
to classical geometers \cite{T}, \cite{G}. The authors of \cite{BN} were
not aware of this.) In the $16$-dimensional
$\F_2$-vector space generated by these lines there is a sub-vector
space of dimension $5$ formed by even sets of lines \cite{N}. All but
the empty set and the sum of all 16 ones consist of eight lines.
So there are 
$$2^5-2 = 30$$
even sets of lines among these $16$ lines.

\vv
$\Lambda_2(8):$ Fix two non-degenerate, disjoint quadrangles
$L_1,...,L_4$ and $L_5,...,L_8$ in space in general position. Both are the
intersections of two smooth quadrics, say
$$L_1+...+L_4=Q_1 \cap Q_2,
\quad L_5+...+L_8=Q_3 \cap Q_4.$$
This shows that the arrangement $L_1,...,L_8$ is the
base locus of the linear system of quartics passing throught it.
For each point on the eight lines there is some quartic
$$(\lambda_1 Q_1 + \lambda_2 Q_2) \cup (\lambda_3 Q_3 + \lambda_4 Q_4)$$
smooth at this point. So again by Bertini we see that the general quartic 
containing these eight lines is smooth.

And of course, the eight lines being two disjoint fibres in the same
elliptic fibration on each of these smooth quartics, are an even
divisor.

\begin{table}[t]
$$
\begin{array}{|c|c|c|} \hline
 & \mbox{arrangement} & k \\ \hline
 & &  \\
\Lambda(6)  & \beginpicture
\put {\line(1,0){40}} [Bl] at 5 -5 
\put {\line(1,0){40}} [Bl] at 55 -5 
\put {\line(2,3){25}} [Bl] at 5 -15
\put {\line(2,-3){25}} [Bl] at 15 22.5 
\put {\line(2,3){25}} [Bl] at 55 -15
\put {\line(2,-3){25}} [Bl] at 65 22.5 
\endpicture & 6 \\
 & & \\ \hline
 & & \\
\Lambda_1(8) & \beginpicture
\put {\line(0,1){40}} [Bl] at 15 -15
\put {\line(0,1){40}} [Bl] at 25 -15
\put {\line(0,1){40}} [Bl] at 35 -15
\put {\line(0,1){40}} [Bl] at 45 -15
\put {\line(0,1){40}} [Bl] at 55 -15
\put {\line(0,1){40}} [Bl] at 65 -15
\put {\line(0,1){40}} [Bl] at 75 -15
\put {\line(0,1){40}} [Bl] at 85 -15
\endpicture & 0 \\
 & & \\
\Lambda_2(8) & \beginpicture
\put {\line(0,1){40}} [Bl] at 10 -15
\put {\line(0,1){40}} [Bl] at 30 -15
\put {\line(1,0){40}} [Bl] at 0 -5
\put {\line(1,0){40}} [Bl] at 0 15
\put {\line(0,1){40}} [Bl] at 70 -15
\put {\line(0,1){40}} [Bl] at 90 -15
\put {\line(1,0){40}} [Bl] at 60 -5
\put {\line(1,0){40}} [Bl] at 60 15
\endpicture & 8 \\
 & & \\
\Lambda_3(8) & \beginpicture
\put {\line(0,1){40}} [Bl] at 35 -15
\put {\line(0,1){40}} [Bl] at 45 -15
\put {\line(0,1){40}} [Bl] at 55 -15
\put {\line(0,1){40}} [Bl] at 65 -15
\put {\line(1,0){40}} [Bl] at 30 -10
\put {\line(1,0){40}} [Bl] at 30 0
\put {\line(1,0){40}} [Bl] at 30 10
\put {\line(1,0){40}} [Bl] at 30 20
\endpicture & 16 \\
 & & \\
\Lambda_4(8) & \beginpicture
\put {\line(5,1){50}} [Bl] at 5 3
\put {\line(5,-1){50}} [Bl] at 5 7
\put {\line(1,1){30}} [Bl] at 25 -5 
\put {\line(1,-1){30}} [Bl] at 25 15 
\put {\line(5,-1){50}} [Bl] at 47 13
\put {\line(5,1){50}} [Bl] at 47 -3
\put {\line(1,-1){30}} [Bl] at 47 25 
\put {\line(1,1){30}} [Bl] at 47  -15 
\endpicture & 16 \\ 
 & & \\ \hline
\end{array} $$
\end{table}

\vv
$\Lambda_3(8):$ 
Let $L_1,...,L_4$ be four distinct lines in one ruling of a smooth
quartic $Q \subset \P_3$ and $L_5,...,L_8$ four distinct lines in the other 
ruling. Grouping them in four pairs like
$$L_1+L_5, \, L_2+L_6, \, L_3+L_7, \, L_4+L_8,$$
each pair generating a plane, we find a quartic (union of the
four planes) cutting out the arrangement $L_1,...,L_8$
on the quadric $Q$. This shows that the arrangement is the
base locus of the linear system of quartics passing through it. For
each point on the arrangement there is a quartic $Q \cup Q'$
smooth at this point. So the general quartic in the system is smooth
again. And the arrangement being cut out on the quartic
by the quadric $Q$ is even.

\vv
$\Lambda_4(8):$ Take two planes $F_1,F_2 \subset \P_3$ meeting in a line
$L$ and choose in each plane a none-degenerate quadrangle, say
$$L_1,L_2,L_3,L_4 \subset F_1, \quad L_5,L_6,L_7,L_8 \subset F_2$$
such that none of the lines $L_i$ coincides with $L$
and such that
$$L_i \mbox{ and }L_{i+4}, \quad i=1,...,4,$$
meet on $L$. Let $E_i$ be the plane spanned by $L_i$ and $L_{i+4}$.
Then $E_1 \cup ... \cup E_4$ is a quartic smooth 
at the four points on $L$ of the arrangement $L_1,...,L_8$. For each other
point on the arrangement there is a quadric surface $Q$ not passing
through this point, such that the quartic $F_1 \cup F_2 \cup Q$
is smooth at this point. Again Bertini shows that the general quartic
through the eight lines is smooth. Beeing cut out by the two
planes $F_1$ and $F_2$ the arrangement is even on each smooth quartic
carrying it.

\VV
\section{Candidates ($n=10$)}

First some terminology:
{\em If the line $L$ in the even arrangement $\Lambda \subset X$ 
meets $l$ other lines from the arrangement, I call the line $L$ 
an $l$-Line. I say the arrangement $\Lambda$ is of the
{\em type} $n_0,n_2,n_4,...$ if $n_0$ is its number of
$0$-lines, $n_2$ is its number of $2$-lines, a.s.o.}

Here I collect the results of section 5.5 below
for arrangements of ten lines: Each even set of ten lines
on a smooth quartic surface is of one of the following eleven types. The
classification in section 5.5 shows that there will not be any other types.
It is not clear to me, however, whether arrangements of the types given
really exist on smooth quartic surfaces. In any case, if such an arrangement
exists, it will be an even set.

\vv
$\Lambda_1(10):$ There are four skew $0$-lines $L_1,L_2,L_3,L_4$ 
and six $2$-lines $L_5,...,L_{10}$ forming
a space hexagon. One line $L$ not belonging to the arrangement 
meets the four $0$-lines, but it does not meet the hexagon.

The six $2$-lines form an elliptic fibre $E$ of degree six. The
divisor $E'=L_1+L_2+L_3+L_4+2L$ forms a $\tilde{D}_4$-fibre in the
same elliptic fibration. By $(\pi 3)$ the divisor $E+E'$ of
degree twelve is even, and the arrangement $E+E'-2L$
then is even too.

\vv
$\Lambda_2(10):$ The union of two disjoint spacial pentagons 
is even by $(\pi 2)$.

\vv
$\Lambda_3(10):$ The arrangement consists of two $0$-lines
$L_1$ and $L_2$, six $2$-lines $L_3,...,L_8$ meeting in pairs
($L_3.L_4=L_5.L_6=L_7.L_8=1$), and two meeting $4$-lines
$L_9,L_{10}$ such that
$$L_3.L_9=L_5.L_9=L_7.L_9=1, \quad
L_4.L_{10}=L_6.L_{10}=L_8.L_{10}=1.$$
There are two skew lines $L,L'$ not belonging to the arrangement, such that
$L$ meets the lines $L_1,L_2,L_3,L_4$, $L'$ meets the lines
$L_1,L_2,L_5,L_6$, with $L,L'$ meeting no other lines from the
arrangement. 

The two disjoint triangles $E=L_3+L_4+L$ and $E'=L_5+L_6+L'$ form two
fibres in the same elliptic fibration of degree three. Elliptic reduction
modulo $E+E'$ reduces the arrangement to $(L_1+L_2+L+L')+
(L_7+L_8+L_9+L_{10})$, the union of two disjoint spacial quadrangles. 
This is an even arrangement of type $\Lambda_2(8)$. So the
original arrangement is even too.

\vv
$\Lambda_4(10):$ The arangement consists of two $0$-lines
$L_1$ and $L_2$, six $2$-lines $L_3,...,L_8$ of which four
meet in pairs ($L_3.L_4=L_5.L_6=1$), and two $4$-lines $L_9,L_{10}$
such that $L_3.L_9=L_5.L_9=L_7.L_9=L_8.L_9=1$ and
$L_4.L_{10}=L_6.L_{10}=L_7.L_{10}=L_8.L_{10}=1$. There are two skew
lines $L,L'$ not belonging to the arrangement such that
$L$ meets $L_1,L_2,L_3,L_4$, the line $L'$ meets $L_1,L_2,L_5,L_6$
and $L,L'$ do not meet any other lines from the arrangement.

Again the two disjoint triangles $E=L_3+L_4+L$ and $E'=L_5+L_6+L'$
are two fibres in the same elliptic fibration. Elliptic reduction modulo
$E+E'$ reduces the arrangement to the even set
$(L_1+L_2+L+L')+(L_7+L_8+L_9+L_{10})$ of type $\Lambda_2(8)$.

\vv
$\Lambda_5(10):$ The arrangement consists of eight $2$-lines
$L_1,...,L_8$ four of which meet in pairs ($L_1.L_2=L_3.L_4=1$)
and two $6$-lines $L_9,L_{10}$ such that
$L_1.L_9=L_3.L_9=L_5.L_9=...=L_8.L_9=1$ and
$L_2.L_{10}=L_4.L_{10}=L_5.L_{10}=...=L_8.L_{10}=1$. There are two lines
$L,L'$ not belonging to the arrangement with $L$ meeting
$L_1,L_3,L_5,...,L_8$ and $L'$ meeting $L_2,L_4,L_5,...,L_8$, while
both $L$ and $L'$ do not meet any other lines from the arrangement.

Again the triangles $E=L_1+L_2+L$ and $E'=L_3+L_4+L'$ are two fibres
in the same elliptic fibration of degree three. Elliptic reduction
modulo $E+E'$ reduces the arrangement to the even
arrangement $L_5+L_6+L_7+L_8+L_9+L_{10}+L+L'$ of type
$\Lambda_3(8)$.

\vv
$\Lambda_6(10):$ There are seven $2$-lines $L_1,...,L_7$ of
which four meet in pairs ($L_1.L_2=L_3.L_4=1$), two skew $4$-lines
$L_8,L_9$ and a $6$-line $L_{10}$ such that the lines
$L_1,...,L_4$ meet the $6$-line $L_{10}$ while the lines
$L_5,L_6,L_7$ meet $L_8$ and $L_9$ with both $4$-lines 
$L_8,L_9$ meeting the $6$-line $L_{10}$.

The $6$-line $L_{10}$ is the intersection of the planes
of the two triangles $L_1+L_2+L_{10}$ and $L_3+L_4+L_{10}$. The
residual lines $L$ (in the plane of $L_1+L_2+L_{10}$) and
$L'$ (in the plane of $L_3+L_4+L_{10}$) therefore are skew, both lines
meeting $L_5,L_5,L_7$. So the two disjoint triangles
$E=L_1+L_2+L$ and $E'=L_3+L_4+L'$ are two fibres in the same
elliptic fibration. Elliptic reduction modulo $E+E'$ reduces
the arrangement to the even set 
$L_5+L_6+L_7+L_{10}+L_8+L_9+L+L'$ of type $\Lambda_3(8)$.

$$ \begin{array}{|c|c|c|c|c|c|c|}
\hline
 & \mbox{arrangement} & k & \qquad &
 & \mbox{arrangement} & k  \\ \hline
 & & & & & & \\
\Lambda_1(10) &
\beginpicture
\put {\line(0,1){40}} [Bl] at 5 -15
\put {\line(3,2){30}} [Bl] at 0 15
\put {\line(3,-2){30}} [Bl] at 0 -5
\put {\line(3,-2){30}} [Bl] at 20 35
\put {\line(3,2){30}} [Bl] at 20 -25
\put {\line(0,1){40}} [Bl] at 45 -15
\put {\line(0,1){40}} [Bl] at 70 -15 
\put {\line(0,1){40}} [Bl] at 80 -15 
\put {\line(0,1){40}} [Bl] at 90 -15
\put {\line(0,1){40}} [Bl] at 100 -15
\setdots<1.5pt>
\plot 65 5  105 5 /
\endpicture & 6 & & \Lambda_6(10) & 
\beginpicture
\put {\line(1,0){110}} [Bl] at 0 20
\put {\line(1,0){30}} [Bl] at 0 10
\put {\line(1,0){30}} [Bl] at 0 0
\put {\line(1,0){30}} [Bl] at 0 -10
\put {\line(0,1){40}} [Bl] at 5 -15 
\put {\line(0,1){40}} [Bl] at 25 -15 
\put {\line(1,-1){20}} [Bl] at 40 25
\put {\line(1,1){20}} [Bl] at 50 5
\put {\line(1,-1){20}} [Bl] at  80 25
\put {\line(1,1){20}} [Bl] at 90 5 
\endpicture & 14  \\
 & & & & & & \\
\Lambda_2(10) &
\beginpicture
\put {\line(1,0){40}} [Bl] at 10 -15
\put {\line(2,-3){20}} [Bl] at 0 10
\put {\line(1,1){35}} [Bl] at 0 0
\put {\line(1,-1){35}} [Bl] at 25 35
\put {\line(2,3){20}} [Bl] at 40 -20
\put {\line(1,0){40}} [Bl] at 95 -15
\put {\line(2,-3){20}} [Bl] at 85 10
\put {\line(1,1){35}} [Bl] at 85 0
\put {\line(1,-1){35}} [Bl] at 110 35
\put {\line(2,3){20}} [Bl] at 125 -20
\endpicture & 10 & & 
\Lambda_7(10) & 
\beginpicture
\put {\line(1,0){70}} [Bl] at -5 35
\put {\line(1,0){50}} [Bl] at 15 15
\put {\line(0,1){70}} [Bl] at 5 -25
\put {\line(0,1){50}} [Bl] at 25 -25
\put {\line(1,0){30}} [Bl] at 0 5
\put {\line(1,0){30}} [Bl] at 0 -5
\put {\line(1,0){30}} [Bl] at 0 -15
\put {\line(0,1){30}} [Bl] at 35 10
\put {\line(0,1){30}} [Bl] at 45 10
\put {\line(0,1){30}} [Bl] at 55 10
\endpicture & 14 \\
 & & & & & & \\
\Lambda_3(10) &
\beginpicture
\put {\line(0,1){50}} [Bl] at 0 -20
\put {\line(0,1){50}} [Bl] at 10 -20
\put {\line(3,1){80}} [Bl] at 20 -20
\put {\line(3,-1){80}} [Bl] at 20 30
\put {\line(1,4){8}} [Bl] at 20 -2 
\put {\line(1,-4){8}} [Bl] at 20 12 
\put {\line(1,4){8}} [Bl] at 40 -2 
\put {\line(1,-4){8}} [Bl] at 40 12 
\put {\line(1,4){8}} [Bl] at 60 -2 
\put {\line(1,-4){8}} [Bl] at 60 12 
\setdots<1.5pt>
\plot -5 10 30 10 /
\plot -5 0 17 0 /
\plot  29 0 50 0 /
\endpicture & 10 & & \Lambda_8(10) &
\beginpicture
\put {\line(1,1){30}} [Bl] at 0 0
\put {\line(1,-1){30}} [Bl] at 0 10
\put {\line(0,1){30}} [Bl] at 10 -10
\put {\line(1,-1){20}} [Bl] at 15 20
\put {\line(1,1){20}} [Bl] at 15 -10
\put {\line(1,0){50}} [Bl] at 20 25
\put {\line(1,0){50}} [Bl] at 20 -15
\put {\line(0,1){50}} [Bl] at 40 -20
\put {\line(0,1){50}} [Bl] at 50 -20
\put {\line(0,1){50}} [Bl] at 60 -20
\endpicture
 & 14\\
 & & & & & & \\
\Lambda_4(10) &
\beginpicture
\put {\line(0,1){50}} [Bl] at 0 -20
\put {\line(0,1){50}} [Bl] at 10 -20
\put {\line(1,0){80}} [Bl] at 20 -10
\put {\line(1,0){80}} [Bl] at 20 20
\put {\line(0,1){50}} [Bl] at 70 -20
\put {\line(0,1){50}} [Bl] at 90 -20
\put {\line(1,4){8}} [Bl] at 20 -2 
\put {\line(1,-4){8}} [Bl] at 20 12 
\put {\line(1,4){8}} [Bl] at 40 -2 
\put {\line(1,-4){8}} [Bl] at 40 12 
\setdots<1.5pt>
\plot -5 10 30 10 /
\plot -5 0 17 0 /
\plot  29 0 50 0 /
\endpicture & 10 & & \Lambda_9(10) & 
\beginpicture 
\put {\line(1,0){60}} [Bl] at 0  30
\put {\line(1,0){60}} [Bl] at 0  -20
\put {\line(0,1){60}} [Bl] at 5 -25
\put {\line(0,1){60}} [Bl] at 55 -25
\put {\line(1,1){20}} [Bl] at 0 15
\put {\line(1,1){20}} [Bl] at 40 -25
\put {\line(2,1){30}} [Bl] at 0 5
\put {\line(1,2){15}} [Bl] at 15 5
\put {\line(1,2){15}} [Bl] at 30 -25
\put {\line(2,1){30}} [Bl] at 30 -10
\endpicture
& 14 \\
 & & & & & & \\
\Lambda_5(10) &
\beginpicture
\put {\line(1,2){18}} [Bl] at 0 0
\put {\line(1,-2){18}} [Bl] at 0 10
\put {\line(1,2){18}} [Bl] at 20 0
\put {\line(1,-2){18}} [Bl] at 20 10
\put {\line(0,1){70}} [Bl] at 50 -30
\put {\line(0,1){70}} [Bl] at 60 -30
\put {\line(0,1){70}} [Bl] at 70 -30
\put {\line(0,1){70}} [Bl] at 80 -30
\put {\line(1,0){85}} [Bl] at 0 -20
\put {\line(1,0){85}} [Bl] at 0 30
\setdots<1.5pt>
\plot 15 2 85 2 /
\plot 0 8  15 8 /
\plot 35 8 85 8 /
\endpicture & 14 & & \Lambda_{10}(10) & 
\beginpicture
\put {\line(0,1){70}} [Bl] at -10 -30
\put {\line(1,0){70}} [Bl] at 0 -25 
\put {\line(1,0){70}} [Bl] at 0 -20 
\put {\line(1,0){70}} [Bl] at 0 -15
\put {\line(0,1){70}} [Bl] at 55 -30
\put {\line(0,1){70}} [Bl] at 60 -30
\put {\line(0,1){70}} [Bl] at 65 -30
\plot 15 -22.5 20 -17.5 /
\put {\line(1,1){27.5}} [Bl] at 25 -12.5
\plot 57.5 20 62.5 25 /
\put {\line(3,1){52.5}} [Bl] at 0 -17.5
\plot 62.5 3.33 67.5 5 /
\put {\line(1,3){17.5}} [Bl] at 40 -12.5
\plot 35 -27.5  36.7 -22.4 /
\endpicture
& 18 \\
 & & & & & & \\ \hline
 & \multicolumn{5}{c|}{ } & \\
\Lambda_{11}(10) & \multicolumn{2}{c}{ 
\beginpicture
\put {\line(0,1){60}} [Bl] at -10 -25 
\put {\line(0,1){6}} [Bl] at 20 -25 
\put {\line(0,1){10}} [Bl] at 20 -17 
\put {\line(0,1){38}} [Bl] at 20 -3 
\put {\line(2,1){38}} [Bl] at 20 -25 
\put {\line(2,1){15}} [Bl] at 65 -2 
\put {\line(2,-1){12}} [Bl] at 20 35 
\plot 35 27.5 38 26 /
\put {\line(2,-1){38}} [Bl] at 42 24 
\put {\line(0,1){6}} [Bl] at 60 -25 
\put {\line(0,1){10}} [Bl] at 60 -17 
\put {\line(0,1){38}} [Bl] at 60 -3 
\put {\line(2,1){38}} [Bl] at 0 5 
\put {\line(2,1){6}} [Bl] at 40 25 
\plot 42 26 45 27.5 /
\put {\line(2,1){12}} [Bl] at 48 29 
\put {\line(2,-1){16}} [Bl] at 0 5 
\put {\line(2,-1){38}} [Bl] at 22 -6 
\put {\line(1,0){8}} [Bl] at 10 -18 
\put {\line(1,0){36}} [Bl] at 22 -18 
\put {\line(1,0){8}} [Bl] at 62 -18 
\put {\line(1,2){6}} [Bl] at 10 -18 
\put {\line(1,2){14}} [Bl] at 18 -2 
\put {\line(1,2){6}} [Bl] at  34 30 
\put {\line(1,-2){6}} [Bl] at  40 42 
\put {\line(1,-2){14}} [Bl] at  48 26 
\put {\line(1,-2){6}} [Bl] at 64 -6 
\endpicture
} & \multicolumn{3}{c|}{\mbox{and two other combinatorical types}} & 18 \\
 & \multicolumn{5}{c|}{ } & \\ \hline
\end{array} 
$$

\vv
$\Lambda_7(10):$ The arrangement consists of six skew $2$-lines
$L_1,...,L_6$ and four $4$-lines $L_7,...,L_{10}$ meting
in pairs ($L_7.L_8=L_9.L_{10}=1$). The $2$-lines $L_1,L_2,L_3$ meet both
$4$-lines $L_7,L_9$ while the $2$-lines $L_4,L_5,L_6$ meet 
both $4$-lines $L_8,L_{10}$. 

The spacial quadrangle $E=L_1+L_2+L_7+L_9$ forms an elliptic fibre. 
The lines $L_8$ and $L_{10}$ are sections for the fibration
$|E|$. The lines $L_4,L_5,L_6$ therefore belong to different fibres 
of $|E|$. Let $E'=L_4+C, \, deg(C)=3$ be the fibre containing
$L_4$. Elliptic reduction modulo $E+E'$ reduces the arrangement to
$(L_3+C)+(L_5+L_6+L_8+L_{10})$. Here $F:=L_5+L_6+L_8+L_{10}$
is an elliptic fibre of degree four. Since $L_3$ is a $2$-section for 
$|E|$, we have $L_3.C=2$. So $F':=L_3+C$ is another elliptic
fibre in $|F|$. Hence $F+F'$ is an even curve, as well as the original
arrangement.

\vv
$\Lambda_8(10):$ The arrangement consists of six $2$-lines
$L_1,...,L_6$ with $L_2.L_3=1$ and four $4$-lines 
$L_7,...,L_{10}$ forming a string ($L_7.L_8=L_8.L_9=L_9.L_{10}=1$).
The line $L_1$ meets $L_8$ and $L_9$, the line $L_2$ meets $L_8$,
the line $L_3$ meets $L_9$ and the three lines $L_4,L_5,L_6$
meet both the lines $L_7$ and $L_{10}$.

The divisor $E:=L_2+L_3+L_8+L_9$ is an elliptic fibre of degree
four with $L_7$ and $L_{10}$ sections for $|E|$. So the lines $L_4,L_5,L_6$
belong to different fibres of $|E|$. Let $E'=L_4+C$ be the fibre in
$|E|$ containing $L_4$. Elliptic reduction modulo $E+E'$
reduces the arrangement to the divisor $F+F'$ with
$F:=L_5+L_6+L_7+L_{10}$ and $F'=L_1+C$ two disjoint elliptic
fibres of degree four. Hence the arrangement is even.
 
\vv
$\Lambda_9(10):$ The arrangement consists of six $2$-lines $L_1,...,L_6$ 
and four $4$-lines $L_7,...,L_{10}$. The $4$-lines form a space quadrangle
($L_7.L_8=L_8.L_9=L_9.L_{10}=L_{10}.L_7=1$). Four $2$-lines meet in pairs
($L_2.L_3=L_5.L_6=1$) while $L_7$ meets $L_1$ and $L_2$,
$L_8$ meets $L_1$ and $L_3$, $L_9$ meets $L_4$ and $L_5$, $L_{10}$
meets $L_4$ and $L_6$.

The divisor $E:=L_2+L_3+L_7+L_8$ is an elliptic
fibre with $L_9,L_{10}$ sections for $|E|$. So $L_4$ and $L_5+L_6$ belong to 
different fibres of $|E|$. Let $E'=L_4+C$ be the fibre containing $L_4$.
Elliptic reduction modulo $E+E'$ reduces the arrangement to $F+F'$ with
$F:=L_5+L_6+L_9+L_{10}$ and $F':=L_1+C$ two disjoint elliptic fibres.
So the original arrangement is even.

\vv
$\Lambda_{10}(10):$ The arrangement contains one $0$-line and nine
$4$-lines. The three $4$-lines $L_2,L_3,L_4$ form a triangle.
Each $4$-line $L_5,L_6,L_7$ meets each of the $4$-lines
$L_8,L_9,L_{10}$. $L_2$ meets $L_5$ and $L_8$, $L_3$ meets 
$L_6$ and $L_9$ while $L_4$ meets $L_7$ and $L_{10}$.

The six lines $L_5,...,L_{10}$ lie on a smooth quadric $Q$.
The residual intersection of $Q$ is a curve $C$ of degree two
with $L_1.C=2$. The curve $C$ does not meet the triangle $L_2+L_3+L_4$.
So quadratic reduction modulo $Q$ reduces the arrangement 
to the union $(L_1+C)+(L_2+L_3+L_4)$ of two disjoint elliptic fibres
of degree three. The original arrangement is even by ???

\vv
$\Lambda_{11}(10):$ Again the arrangement contains one
$0$-line $L_1$ and nine $4$-lines $L_2,...,L_{10}$.
These nine $4$-lines form three triangles $L_2+L_3+L_4$,
$L_5+L_6+L_7$ and $L_8+L_9+L_{10}$. Each line from one triangle meets
precisely one line from each of the two other triangles.
We can reorder these lines such that
$$L_2.L_5=L_2.L_8=1, \, L_3.L_6=L_3.L_9=1, \, L_4.L_7=L_4.L_{10}=1.$$
Then there are three essentially different combinatoriv\cal possibilities
for the intersection pattern of the lines $L_5,L_6,L_7$ with the lines
$L_8,L_9,L_{10}$. 

Let $L$ resp. $L'$ be the residual intersections of the planes of the
triangles $L_2+L_3+L_4$, resp. $L_5+L_6+L_7$. Both the lines $L,L'$ are
skew with the third triangle while both of them meet $L_1$. 
Quadratic reduction modulo the two planes 
leads to the even set $(L_1+L+L')+(L_8+L_9+L_{10})$ of type
$\Lambda(6)$. Hence the original arrangement is even.

\VV
\section{Classification}

Here I want to classify even arrangements $\Lambda = \sum_1^n L_i$ of 
$n\leq 8$ lines
$L_i \subset X 
\subset \P_3$ on a smooth quartic surface $X$.

\vV
\subsection{Two lines}

For $n=2$ there are only the possibilities $k=0,1$ violating
the modulo-$4$-condition $(\lambda 6)$.

\vV
\subsection{Four lines}

For $n=4$ we cannot have four skew lines, because this violates 
$(\pi 1)$. If they are not skew, 
by ($\lambda 6$) we have $k=4$, and each line must meet
two other lines. All four of them form a space quadrangle. This is an
elliptic fibre and cannot be even by $(\pi 2)$.  

\vV
\subsection{Six lines}

For $n=6$ we have $k \leq 12$, hence $k=2,6$ or $10$. In all cases
$\LL$ is effective by Riemann-Roch.

{\bf k=2:} There are just two $2$-lines meeting in at most
one point, impossible.

\vv
{\bf k=6:} We have $\Lambda^2=0$. So $\LL$ is effective of degree three
with $\LL^2=0$.  If it were not reduced, it would contain a
multiple line $L$. If $\LL \sim 3L$, then $\LL^2 = -18$, impossible.
If $\LL=2L+L'$, the $\LL^2=-6$ or $=-10$, impossible too.
So $\LL$ is represented by a reduced divisor $E$. This divisor must
be connected, because $\LL^2 < 0$ otherwise. Since $deg(E)=3$,
the linear system $|E|$ cannot have a fixed component.
So $E$ is an elliptic
fibre with $\Lambda \sim 2E$. This implies that
$\Lambda$ consists of two elliptic fibres of degree three. 
This is the type $\Lambda(6)$.   

\vv
{\bf k=10:} We have
$$\chi(\OO_X(2) \otimes \OO_X(-\Lambda))>0.$$
Since
$$deg(\OO_X(\Lambda) \otimes \OO_X(-2))=-2<0,$$
necessarily
$$h^2(\OO_X(2) \otimes \OO_X(-\Lambda))=
h^0(\OO_X(\Lambda) \otimes \OO_X(-2))=0.$$
So the arrangement $\Lambda$ lies on a quadric surface $S$.
If $S$ is smooth, the arrangement can consist of six
skew lines ($k=0$), or four lines in one ruling and two from the other
one ($k=8$). In both cases $k \not= 10$. If $S$ is not smooth,
but breaks up into two planes, each plane contains at most four
lines and there are the following two cases:  

- Either one plane contains two lines only, the other one four. On these 
four lines there are altogether four points to be met by the two lines
from the other plane.

- Or both planes would contain three lines.  $k=10$ implies that each of them 
meets a line from the other plane, impossible, because then each line
meets exactly three other ones.

\vV
\subsection{Eight lines}

For $n=8$ we have $k \leq 24$, so 
$$k=0,4,8,12,16,20 \mbox{ or } 24.$$
Here $k=0$ is possible: Type $\Lambda_1(8)$. In all other cases 
$\LL$ is effective by Riemann Roch, with $deg(\LL)=4$
and $\LL^2=\frac{k}{2}-4$.

{\bf k=4:} There cannot be a $\geq 4$-line, because it would meet
four none-$0$-lines causing $k>4$. So there will be four
$2$-lines $L_1,...,L_4$ forming a space quadrangle and 
four $0$-lines $L_5,...,L_8$. For each of them $\LL.L_i=-1$,
hence $\LL$ splits off the $0$-line $L_i$. This implies
$\LL\sim L_5+...+L_8$ and $\LL^2=-8$, a contradiction.

\vv
{\bf k=8:} Now $\LL$ is effective with $\LL^2=0$
and $deg(\LL)=4$. By Riemann-Roch $h^0(\LL) \geq 2$.
If $|\LL|$ has no fixed component, it is an elliptic fibration
of degree four. And $\Lambda \sim 2E$ will consist of two fibres 
in $E$. Both fibres are space quadrangles and we obtain case
$\Lambda_2(8)$.

If however $|\LL|$ has a fixed component $L$, then necessarily
$deg(\LL-L) \geq 3$, so $L$ is a line and $\LL \sim L+E$
with $|E|$ an elliptic fibration of degree three. From
$$0=\LL^2=L^2+2L.E+E^2=2L.E-2$$
we deduce $E.L=1$. Hence
$$\LL.L =(L+E).L=-1 \quad \mbox{and} \quad \Lambda.L=-2.$$
Therefore $L=L_1$ is a $0$-line in $\Lambda$. Further
$$E.(\Lambda-L_1)=E.(2E+L_1)=1$$
shows that there is exactly one line $L_8 \subset \Lambda$ with
$$E.L_8=1, \quad \LL.L_8=(E+L_1).L_8=1, \quad 
\Lambda.L_8=2.$$
So $\Lambda$ contains exactly one $4$-line $L_8$ and six $2$-lines
$L_2,...,L_7$. The line $L_8$ is a section for $|E|$ and the lines $L_2,...,L_7$
are fibre components of $|E|$. A fibre cannot consist of three lines
$L_i, i=2,...,7$, because it also meets $L_8$. It cannot contain
two lines $L_i,L_j, \, i,j=2,...,7$ from the arrangement $\Lambda$,
because one of them would meet only one line from $\Lambda$. And a
fibre cannot contain one line $L_i, i=2,...,7,$ only, because this $2$-line
would meet at most the line $L_8$ from $\Lambda$. Contradiction!

\vv
{\bf k=12:} Now $\LL$ is effective with $deg(\LL)=4$
and $\LL^2=2$. If there is a $0$-line $L_1$ belonging to $\Lambda$,
with $\LL.L_1=-1$ showing that $\LL$ splits off this line.
So $\LL \sim L_1+C$ with $C$ effective of degree three and
$\LL.L_1=-2+L_1.C=-1$. This implies $L_1.C=1$. But then
we arrive at the contradiction
$$\LL^2 = -2+2L_1.C+C^2=C^2 \leq 0.$$

Assume next that $\Lambda$ contains two $2$-lines
$L_1,L_2$. Then
$\Lambda.L_i=0$ implies $\LL.L_i=0$. From Riemann Roch we find
$h^0(\LL)\geq 3$, so there is a divisor $D \sim \LL$ meeting $L_1$
and $L_2$. 
By $D.L_i=0$ this divisor splits off $L_1$ and $L_2$, say $D=L_1+L_2+C$
with $deg(C)=2$. We have
$$0=D.L_1=-2+L_1.L_2+L_1.C.$$
If $L_1.L_2=1$, then $L_1.C=L_2.C=1$ and
$$D^2=(L_1+L_2+C)^2\leq -6+6=0,$$
in conflict with $\LL^2=2$. If $L_1.L_2=0$, then $L_1.C=L_2.C=2$ and
$$D^2=(L_1+L_2+C)^2=4+C^2=2$$ 
only if $C^2=-2$. This shows that $C$ is planar, i.e.\ a smooth conic or two 
intersecting lines. Let $P \subset \P_3$ be the plane containing $C$.
Then 
$$L_i.P \geq L_i.C=2, \, i=1,2,$$
shows that $L_i \subset P$ too. So $D \sim \LL$ is a plane section
with $D^2=4$, again a contradiction.

So the type of $\Lambda$ is
$$(0,n_2,n_4,n_6)$$
with $n_2 \leq 1$, hence $n_4+n_6 \geq 7$. But then
$\sum k_i \geq 28$ and $k \geq 14$, a contradiction! 

\vv
{\bf k $\geq$ 16:} Riemann-Roch as usual shows that the whole
arrangement $\Lambda$ lies on a quadric $Q$. If $Q$ is smooth,
we have four lines in each of the rulings, type $\Lambda_3(8)$.
If $Q$ breaks up into two planes, we have four lines in
each plane, type $\Lambda_4(8)$.

\vV
\subsection{Ten lines}

Now $k \leq 40$, hence
$$k=2,6,10,...,38.$$
However
$$\chi(\OO_X(2)-\Lambda)=2+\frac{1}{2}(16-40+2k-20)=
k-20 > 0,$$
if $k>20$. Since $deg(\Lambda-\OO_X(2))=2$, if $C=\Lambda-\OO_X(2)$
is effective, this divisor $C$ can consist of two lines,
skew or meeting, or a smooth conic. In all three cases
$$h^0(C)=h^2(\OO_X(2)-\Lambda)=1.$$
So for $k \geq 22$ we have the contradiction
$$h^0(\OO_X(2)-\Lambda)>0.$$
This shows that there are in fact only the possibilities
$$k=2,6,10,14,18.$$

The case  k=2 is impossible as usual. In all other cases
$k>n-8$ and $\LL$ is effective by section 2. In fact
$$deg(\LL)=5, \quad \LL^2 =\frac{1}{2}(k-10).$$

{\bf k=6:} The arrangement $\Lambda$ cannot contain a $6$-line, because the
six lines meeting it would increase $k$ at least by three to 
give $k \geq 9$. So we have
$$2 \cdot n_2 + 4 \cdot n_4=12.$$
If $n_4=0$, there are four $0$-lines $L_0,...,L_4$ and six $2$-lines
$L_5,...,L_{10}$. These six $2$-lines cannot form two disjoint 
triangles, because this would be an even configuration $\Lambda(6)$,
and subtracting it from $\Lambda$, 
we would find the contradiction that the divisor
$L_1+...+L_4$ is even. So $E=L_5+...+L_{10}$ is a spacial hexagon.
Now $\LL.L_i=-1$ for the four $0$-lines $L_i,\, i=1,...,4$. So
$\LL$ splits them off, say
$$\LL=L_1+...+L_4+L$$
with a fifth line $L \subset X$. Since $\LL.L_i=0$ for 
$i=5,...,10$, this line $L$ does not meet the hexagon $E$.
This is type $\Lambda_1(10)$.

If $n_4=1$, there are five $0$-lines $L_1,...,L_5$, four
$2$-lines $L_6,...,L_9$ and the $4$-line $L_{10}$. Now
$\LL=L_1+...+L_5$ with $\LL.L_{10}=1$, a contradiction.

If $n_4 \geq 2$, we would have at least six $0$-lines contained in $\LL$,
a contradiction.

\vv
{\bf k=10:} Now $\LL^2=0$. If $n_0 \geq 3$, then
$$\LL=L_1+L_2+L_3+C$$
with $0$-lines $L_1,L_2,L_3$ and a curve $C \subset X$ of
degree two. Since $\LL.L_i=-1$ we find $L_i.C=1$ for $i=1,2,3$.
But since $C^2 \leq -2$, this leads to the contradiction
$\LL^2 \leq -2$. Hence $n_0 \leq 2$.

If $n_0=0$, the configuration $\Lambda$ consists of
cycles of $2$-lines. It cannot be a single cycle only by ($\pi 2$),
so it must consist of two cycles of the same length five,
type $\Lambda_2(10)$.

If $n_0=1$, we have $n_2=8$ and $n_4=1$. The $4$-line $L_{10}$ meets
four $2$-lines. Each of them belongs to a string 
$L_i,L_{i+1},...,L_{i+i'}$ of $2$-lines with 
$$L_i.L_{10}=1, \quad L_j.L_{j+1}=1 \mbox{ for }j=i,...,i+i'-1,
\quad L_{i+i'}.L_{10}=1.$$
Since $h^0(\LL)\geq2$ and $\LL.L_i=0$, there is a divisor 
$D \sim \LL$ splitting off $L_1$ and $L_i$. Then it splits off 
$L_{i+1},...,L_{i+i'}$ too and from $\LL.L_{10}=1$ we conclude
$L_{10} \subset D$. Then $D$ would also split off the other two 
$2$-lines meeting $L_{10}$,
and $\LL$ would be linearly equivalent to the sum 
of at least six lines, contradiction! 

If $n_0=2$, we have 
$$n_2=7, \, n_4=0, \, n_6=1, \quad \mbox{ or } \quad
n_2=6, \, n_4=2, \, n_6=0.$$
Now $\LL$ splits off the two $0$-lines, say
$$\LL \sim L_1+L_2+C$$ 
with $C$ a curve of degree three satisfying $L_1.C=L_2.C=1$. But then
$$0=\LL^2=L_1^2+L_2^2+C^2+2(L_1+L_2).C=C^2$$
implies that $|C|$ is an elliptic fibration of degree three. If there
were a $6$-line $L_{10}$, it would be the line residual to the fibration.
Then $C.L_{10}=3$, a contradiction. This shows $n_2=6, n_4=2$.

The two $0$-lines $L_1$ and $L_2$ are sections for
the elliptic fibration $|C|$, while each $2$-line will be 
a fibre component. There cannot be three $2$-lines in a fibre
(at least one of them would meet $L_1$ or $L_2$), so each $2$-line
meets at most one other $2$-line. And if two $2$-lines meet, 
then they cannot both meet the same
$4$-line $L_{10}$, because $\LL.L_{10}=1$ shows that $L_{10}$ 
is a section for $|C|$.
There are the 
following two possibilities:

$$
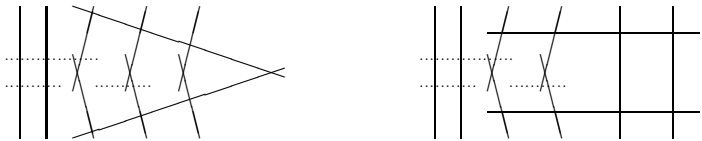
\beginpicture
\put {\line(0,1){50}} [Bl] at 0 0
\put {\line(0,1){50}} [Bl] at 10 0
\put {\line(3,1){80}} [Bl] at 20 0
\put {\line(3,-1){80}} [Bl] at 20 50
\put {\line(1,4){8}} [Bl] at 20 18 
\put {\line(1,-4){8}} [Bl] at 20 32 
\put {\line(1,4){8}} [Bl] at 40 18 
\put {\line(1,-4){8}} [Bl] at 40 32 
\put {\line(1,4){8}} [Bl] at 60 18 
\put {\line(1,-4){8}} [Bl] at 60 32 
\setdots<2pt>
\plot -5 30 30 30 /
\plot -5 20 17 20 /
\plot  29 20 50 20 /
\endpicture 
\hspace{2cm}
\beginpicture
\put {\line(0,1){50}} [Bl] at 0 0
\put {\line(0,1){50}} [Bl] at 10 0
\put {\line(1,0){80}} [Bl] at 20 10
\put {\line(1,0){80}} [Bl] at 20 40
\put {\line(0,1){50}} [Bl] at 70 0
\put {\line(0,1){50}} [Bl] at 90 0
\put {\line(1,4){8}} [Bl] at 20 18 
\put {\line(1,-4){8}} [Bl] at 20 32 
\put {\line(1,4){8}} [Bl] at 40 18 
\put {\line(1,-4){8}} [Bl] at 40 32 
\setdots<2pt>
\plot -5 30 30 30 /
\plot -5 20 17 20 /
\plot  29 20 50 20 /
\endpicture$$

\noindent
Each pair of meeting $2$-lines can be completed by an additional line
(dotted) to a fibre of $|C|$. 
These are the types $\Lambda_3(10)$ and $\Lambda_4(10)$.

\vv
{\bf k=14:} Now $\LL^2=2$ and by Riemann-Roch $h^0(\LL) \geq 3$. First
of all we observe that there cannot be any $0$-line in the
arrangement $\Lambda$. If there were one, say $L_1$, then 
$\LL \sim L_1+C$ with a curve $C$ of degree four. From
$\LL.L_1=-1$ we conclude
$$C.L_1=(\LL-L_1).L_1=1, \quad C^2=(\LL-L_1)^2=2.$$
As $deg(C)=4$ we have $(C-\OO(1))^2=6-8=-2$, and
either $C-\OO(1)$ or $\OO(1)-C$ is effective. In both
these cases $C \sim \OO(1)$ in conflict with $C^2=2$. 

Again
it is impossible, that a $2$-line, say $L_2$, meets two other
ones, say $L_1$ and $L_3$. Because then there is a divisor
$$D= L_1+L_2+L_3+C \sim \LL$$
with a curve $C$ of degree two. There are two cases:

Either $L_1$ and $L_3$ also meet. Then from $\LL.L_i=0$
we conclude, that $C.L_i=0$. But $C^2 < 0$ then would imply
the contradiction $\LL^2 < 0$. 

Or $L_1$ and $L_3$ don't meet. Then
$$C.L_1=C.L_3=1, \quad C.L_2=0.$$
Again $C^2<0$ would imply
$$\LL^2=(L_1+L_2+L_3)^2+2(L_1+L_3).C + C^2=2+C^2 \leq 0.$$

The two equations
\begin{eqnarray*}
2n_2+4n_4+6n_6+8n_8 &=& 28 \\
2n_2+2n_4+2n_6+2n_8 &=& 20 
\end{eqnarray*}
show
$$2n_4+4n_6+6n_8=8.$$
There are the following four possibilities:
$$ \begin{array}{c|cccc}
\mbox{case} & n_2 & n_4 & n_6 & n_8 \\  \hline
a)          & 8   & 1   & 0   & 1   \\
b)          & 8   & 0   & 2   & 0   \\
c)          & 7   & 2   & 1   & 0   \\
d)          & 6   & 4   & 0   & 0   \\
\end{array}$$ 

{\bf Observation:} Each $2$-line in the arrangement meets at most one other
$2$-line. If there were two such pairs of meeting $2$-lines,
say $L_1.L_2=L_3.L_4=1$, then by $h^0(\LL) \geq 3$ there would
be a divisor $L_1+L_2+L_3+L_4+L \sim \LL$ with a fifth line $L \subset X$.
From $\LL.L_i=0$ follows
$$\LL.L=\LL.(\LL-(L_1+...+L_4))=2, \quad
L_i.L=L_i.(\LL-(L_1+...+L_4))=-1-L_i^2=1.$$
So either $L$ belongs to $\Lambda$, being a $6$-line, or
$L$ does not meet any line $L_i, \, i>4$. Then each $4$-line meets one of the 
lines $L_i, \, i \leq 4$, each $6$-line two of them and each $8$-line
three.

Case a): If the $4$-line $L_9$ meets the $8$-line $L_{10}$,
then on $L_9+L_{10}$ there are ten points to be met
by eight $2$-lines. So there will be one $2$-line,
say $L_1$, which meets both the lines $L_9$ and $L_{10}$.
On $L_9+L_{10}$ there remain eight points to be met by
seven $2$-lines. So another $2$-line, say $L_2$ meets both
$L_9$ and $L_{10}$. Then $L_1$ and $L_2$, being coplanar,
will meet, contradiction!

If $L_9$ is skew with $L_{10}$, then only four points on $L_{10}$
can be connected with the four points on $L_9$ by $2$-lines.
By the observation above among them there is a string of two. Let 
$L_1,...,L_4$ be the four $2$-lines meeting $L_9$ and $L_4+L_5$
be the string. By $h^0(\LL) \geq 3$ and $\LL.L_9=1$ there is a
divisor $D \sim \LL$ splitting off $L_9$. Then it splits off
$L_1,...,L_5$ too and consists of six lines, contradiction.

Case b): If the two $6$-lines, say $L_9$ and $L_{10}$ meet, 
there are ten points on $L_9+L_{10}$ to be met by eight
$2$-lines. One $2$-line must lie in the plane
spanned by $L_9$ and $L_{10}$. Then there remain on $L_9+L_{10}$
eight points to be met by seven $2$-lines. A second $2$-line must
lie in this plane, contradiction!

If $L_9$ and $L_{10}$ are skew, their twelve points must
be connected by eight $2$-lines or strings of those. There
will be exactly two strings $L_1+L_2$ and $L_3+L_4$ of $2$-lines. By the 
observation above there will be a line $L \subset X$ meeting
these four lines. If $L$ is one of the
$6$-lines, say $L=L_{10}$, then $L_{10}$ meets all four lines
$L_i, \, i=1,...,4$. Only two more points
on $L_{10}$ can be connected with the six points on $L_9$
by $2$-lines, contradiction. So $L$ does not belong to
$\Lambda$ and does not meet any line $L_i,\, i>4$. This implies
that not both the lines in one string meet the same
$6$-line.

Let $L'\subset X$ be the residual line in the plane
of $L_1,L_2,L$ and $L'' \subset X$ the residual line in the plane
of $L_3,L_4,L$. Both the lines $L'$ and $L''$ are skew,
skew with $L_9$ and $L_{10}$, meeting $L_5,...,L_8$, 
type $\Lambda_5(10)$. 

$$
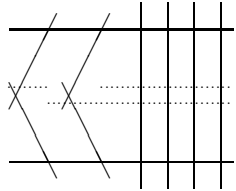
\beginpicture
\put {\line(1,2){18}} [Bl] at 0 30
\put {\line(1,-2){18}} [Bl] at 0 40
\put {\line(1,2){18}} [Bl] at 20 30
\put {\line(1,-2){18}} [Bl] at 20 40
\put {\line(0,1){70}} [Bl] at 50 0
\put {\line(0,1){70}} [Bl] at 60 0
\put {\line(0,1){70}} [Bl] at 70 0
\put {\line(0,1){70}} [Bl] at 80 0
\put {\line(1,0){85}} [Bl] at 0 10
\put {\line(1,0){85}} [Bl] at 0 60
\setdots<2pt>
\plot 15 32 85 32 /
\plot 0 38  15 38 /
\plot 35 38 85 38 /
\endpicture
$$

Case c): Let $L_8$ and $L_9$ be the two $4$-lines and $L_{10}$
the $6$-line. On them there are altogether $14$ points, where they
can meet other lines. If a $4$-line, say $L_8$ is skew with
$L_9$ and $L_{10}$, then it meets four $2$-lines, say $L_1,...,L_4$.
By $h^0(\LL)\geq 3$ there is a divisor $D=L_1+...+L_4+L_8 \sim \LL$.
Since $D^2=2$ the four $2$-lines $L_1,...,L_4$ meet in pairs.
But this causes the contradiction $\Lambda.L_9=\Lambda.L_{10}=0$.
So each $4$-line meets the other $4$-line or $L_{10}$.

Assume that $L_8$ meets $L_9$. Then $L_{10}$ cannot lie in their
common plane, because this plane would have to contain a
$2$-line too, contradiction.

Assume that $L_{10}$ meets meets one of the $4$-lines, say
$L_9.L_{10}=1$. No $2$-line, or a string of those, can connect
points on $L_8$ with points on $L_9$, because by $h^0(\LL)=3$
there would be a divisor $D \sim \LL$ splitting off $L_8$, as well as the
at least three $2$-lines meeting it, then $L_9$ too with at least
another $2$-line. These are together at least six lines,
contradiction. Not all three points on $L_8$ can be
connected with points on $L_{10}$, because then there would be a divisor
$D \sim \LL$ splitting off $L_8+L_{10}$ and all the lines
meeting $L_{10}$. These are too many again. There remains only the
possibility, that one string of $2$-lines connects
two points on $L_8$, one $2$-line connects a point on $L_8$ with a point
on $L_{10}$, the second string connects two points on $L_{10}$, while
two $2$-lines connect points on $L_{10}$ with points on $L_9$.
So the plane of $L_9$ and $L_{10}$ contains two $2$-lines,
contradiction.

So if $L_8$ meets $L_9$, the line $L_{10}$ will be skew
with both of them. If a point on $L_8$ is connected with a point
of $L_9$ by $2$-lines, there will be a divisor $D \sim \LL$
splitting off $L_8+L_9$ as well as the at least five $2$-lines meeting them,
impossible.

This shows $L_8.L_9=0$, hence $L_8.L_{10}=L_9.L_{10}=1$. No two points on $L_8$
(or $L_9$) can be connected with two points on $L_{10}$,
because then a divisor $D \sim \LL$ would split off $L_8+L_{10}$ and
the four $2$-lines meeting $L_{10}$, impossible. This implies that 
two points on $L_{10}$ will be connected by a string $L_1+L_2$
of $2$-lines. Let $L_3,L_4$ be the other two $2$-lines meeting
$L_{10}$. By $h^0{\LL} \geq 3$ some divisor $D \sim \LL$ will split off
$L_1+L_2$ and $L_3$. Then $D^2=2$ implies $L_3.L_4=1$ and
$D=L_1+L_2+L_3+L_4+L_{10}$. The remaining three points on $L_8$ and on $L_9$ 
then are joined by the remaining three $2$-lines $L_5,L_6,L_7$,
type $\Lambda_6(10)$.

Case d): Let $L_1,...,L_6$ be the $2$-lines and $L_7,...,L_{10}$
the $4$-lines. If a $4$-line doesn't meet another $4$-line, it
will meet four $2$-lines, and as above we find a divisor
$D \sim \LL$ with $D^2\leq 0$, impossible, or $D.L_i=0$ for the
other three $4$-lines $L_i$, impossible too. 

If the four $4$-lines meet in pairs, say $L_7.L_8=L_9.L_{10}=1$, 
otherwise skew, then there cannot be a $2$-line, say $L_1$, meeting both
$L_7$ and $L_8$. Because then there will be an effective
divisor splitting off $L_1+L_7+L_8$, as well as the 
at least four other
$2$-lines meeting $L_8+L_9$, contradiction. So the remaining six points
on $L_7+L_8$ will be joined with the remaining six points on $L_9+L_{10}$ by
$2$-lines. Let $L_1,L_2,L_3$ be the $2$-lines meeting
$L_7$. At least two of them, say $L_1,L_2$ will meet the same other $4$-line,
say $L_9$. Then there is an effective divisor $D \sim \LL$
splitting off $L_1+L_2+L_7+L_9$. If $L_3$ doesn't meet $L_9$, then
$D$ splits off $L_8$ too, as well as at least two other $2$-lines,
contradiction. Hence $L_1,L_2,L_3$ join points on $L_7$ with points on
$L_9$ while $L_4,L_5,L_6$ join points on $L_8$ with points on $L_{10}$,
type $\Lambda_7(10)$.

We are left with the cases where $d$, the number of intersection points
of $4$-lines is at least $=3$. Then there are precisely $d-2$ strings
of meeting $2$-lines. By the observation above, $d \leq 4$.
We consider the possibilities:

d=3: The four $4$-lines form one string, say 
$L_7.L_8=L_8.L_9=L_9.L_{10}=1$. If a $2$-line $L_1$ connects two points
on $L_i$ and $L_{i+1}, \, i \geq 7$, then there is a divisor $\sim \LL$
splitting off $L_1+L_i+L_{i+1}$ as well as all the other 
$2$-lines meeting $L_i+L_{i+1}$. This is possible only if
$L_1$ meets $L_8$ and $L_9$. Let $L_2$ be the second $2$-line meeting
$L_7$ and $L_3$ the second $2$-line meeting $L_8$. Then 
$\LL \sim L_1+L_2+L_3+L_8+L_9$. We get the contradiction
$\LL^2=0$ unless $L_2.L_3=1$. This is type $\Lambda_8(10)$.  

d=4: Assume first that the four $4$-lines form a spacial
quadrangle. If a $2$-line $L_1$ connects two points on meeting
$4$-lines, say $L_7$ and $L_8$, we get the same contradiction as above, unless
another string $L_2+L_3$ connects the remaining points on $L_1+L_2$.
Then the remaining points on $L_9+L_{10}$ must be connected 
by a line $L_4$ and a string $L_5+L_6$, type $\Lambda_9(10)$. 
  
In all other cases each $4$-line $L_i$ will be connected to 
the opposite one by two $2$-lines and one string of those. Then
$\LL$ would split off all ten lines, contradiction.

If the $4$-lines do not form a spacial quadrangle, there will be one
of them, say $L_7$ meting just one other $4$-line, say $L_8$, while
$L_8+L_9+L_{10}$ form a triangle. On $L_8$ there remains a fourth point
to be connected by one or two $2$-lines with a point on 
another $4$-line $L_i$. As $L_8$ meets all other $4$-lines,
there will be a divisor $\sim \LL$ splitting off $L_8,L_i$
and all the $2$-lines meeting this pair. If $L_i$ belongs to
the triangle, this divisor will split off the whole triangle, and then all
the lines, too much. So $L_i=L_7$. There remain two points on $L_7$.
If they are connected by a string of $2$-lines, there
remain two points on each of the lines $L_9,L_{10}$. Only two of them can
be connected by the second string, two must be connected by one line.
This $2$-line then lies in the plane of the triangle, contradiction.
The remaining points on $L_7$ therefore are connected to points
on $L_9+L_{10}$. The divisor will split off at least one of them, and then
all the lines, contradiction!

\vv
{\bf k=18:} Now $\LL^2=4$ and $h^0(\LL) \geq 4$ by Riemann-Roch.

Let's first consider the (strange) case that $\Lambda$ contains
a $0$-line $L_1$. Then each divisor $D \sim \LL$ splits off
$L_1$, i.e. $D=L_1+C$ with $deg(C)=4$. From $\LL.L_1=-1$
we conclude $C.L_1=1$
and $C^2=(D-L_1)^2=4$. The linear
system $|C|$ cannot have a fixed component, because then 
$h^0(\LL)=h^0(C)\leq 2$. So by Bertini the generic $C \in |C|$
is irreducible with arithmetical genus $=3$. This is impossible, if
$C$ is a genuine space curve. So $C$ is a plane section of the
surface and $\LL \sim L_1+\OO(1)$. 

For all the other nine lines $L_i, i=2,...,10,$ in $\Lambda$ this implies
$\Lambda.L_i=2$, i.e. they are $4$-lines. We immediately observe that
there cannot be a planar quadrangle of $4$-lines, because this would
meet the line $L_1$. Next we claim that each $4$-line belongs to at least
one triangle of $4$-lines. Indeed, if 
the $4$-line $L_2$ would meet four $4$-lines $L_3,...,L_6$
which don't intersect, on $L_2+L_3+...+L_6$ there would be
16 points of intersection. There remain two points, where 
$L_7,...,L_{10}$ can meet. So either these lines meet in pairs,
say $L_7.L_8=L_9.L_{10}=1$ or there is a string of three lines,
say $L_7.L_8=L_8.L_9=1$. In the first case, both the lines $L_7$ and $L_8$ 
would meet three lines from $L_3,...,L_6$. Then $L_7$ and $L_8$ would form
a planar quadrangle with two of the lines from $L_3,...,L_6$, contradiction.
In the second case 
both the lines $L_7$ and $L_9$ would meet three lines from $L_3,...,L_6$
forming a planar quadrangle with two of them, again a contradiction.

So we may assume that the lines $L_2,L_3,L_4$ form a triangle. Let
$L \subset X$ be the residual line of their common plane. On each of
these three lines there are two more points of intersection with one of the
six lines $L_5,...,L_{10}$. This implies that each of the lines
$L_i, \, i=5,...,10$ meets meets three other ones from those lines. we
distinguish two cases:

Either there is no triangle among $L_5,...,L_{10}$. Then they can be 
grouped into two triplets, say $L_5,L_6,L_7$ and $L_8,L_9,L_{10}$ 
such that each line from the first triplet meets each line from the
second one, type $\Lambda_{10}(10)$. 

Or there is a second triangle, say $L_5,L_6,L_7$. Each line $L_i, \,
i=5,6,7$ meets exactly one line from the triangle $L_1,L_2,L_3$.
The lines $L_8,L_9,L_{10}$ therefore form a triangle too.
This is case $\Lambda_{11}(10)$ containing three combinatorically
different subcases. 

From now on we may assume that $\Lambda$ does not contain any
$0$-line. But $\Lambda$ will not contain any $2$-line either: If
there is such a $2$-line $L_1$ then $\LL \sim L_1+C$ with $C$
effective of degree four. From $\LL.L_1=0$ we conclude
$$C^2=(\LL-L_1)^2=2.$$
Hence $(C-\OO(1))^2=-2$ and either $C-\OO(1)$ or $\OO(1)-C$ is
effective. In both cases $C\sim \OO(1)$ in conflict with $C^2=2$.

So we have the two equations
\begin{eqnarray*}
4n_4+6n_6+8n_8 &=& 36 \\
2n_4+2n_6+2n_8 &=& 20
\end{eqnarray*}
leading to the contradiction
$$n_6+2n_8=-2.$$
There are no arrangements $\Lambda$ of this type.

\VV
\section{Chern numbers}

It is tempting to compute the Chern numbers of the double covering
surface and to apply the known conditions for these numbers.
I shall do this in this section, although the result is disappointing.

For simplicity, let me assume that {\em no three (or four) lines in
the arrangement $\Lambda$ are concurrent}. Then consider the
succession of maps:
$$ X \leftarrow \tilde{X} \leftarrow \tilde{Y} \rightarrow Y.$$
The maps are:

$\tilde{X} \to X$ is the blow up of the quartic surface $X$
in the $k$ points, where lines from the arrangement intersect.
Over each point there is introduced some $(-1)$-curve $E_j$. 
Any $l$-line $L \subset \Lambda \subset X$ corresponds to a rational
curve $M \subset \tilde{X}$ with self-intersection 
$-2-l$. The canonical divisor of $\tilde{X}$ is
$K_{\tilde{X}}=\sum E_j$. The Chern numbers of $\tilde{X}$ are
$$c_1^2(\tilde{X})=\sum E_j^2=-k, \quad 
c_2(\tilde{X})=c_2(X)+k = 24+k.$$

$\tilde{Y} \to \tilde{X}$ is the double cover branched over 
the $n$ rational curves $M_i$. The $(-1)$-curves $E_j$ 
correspond to $(-2)$-curves $\tilde{E}_j \subset \tilde{Y}$. 
The rational curves $M_i$ correspond to rational curves
$\tilde{M}_i \subset \tilde{Y}$ of self-intersection
$$\frac{1}{2}(-2-l_i)=-1-\frac{l_i}{2}.$$
The canonical divisor of $\tilde{Y}$ is
$$K_{\tilde{Y}} = \sum \tilde{E}_j+\sum \tilde{M}_i.$$
The Chern numbers of $\tilde{Y}$ therefore are
$$c_1^2(\tilde{Y})=(\sum\tilde{E}_j+\sum\tilde{M}_i)^2=
-2k+2\sum l_i -\sum(1+\frac{l_i}{2})=k-n, \quad
c_2(\tilde{Y})=2 c_2(\tilde{X})-2n = 48+2(k-n).$$

If the arrangement $\Lambda$ contains $0$-lines, the surface $\tilde{Y}$
is not minimal. One obtains a minimal surface $Y$ by blowing down
the $(-1)$-curves $\tilde{M}_i \subset \tilde{Y}$ 
corresponding to $0$-lines $L_i \subset X$
via the third map $\tilde{Y} \to Y$. Denote by 
$\tilde{L}_i \subset Y$ the images of the other curves
$\tilde{M}_i$ and by $F_j \subset Y$ the images
of the $(-2)$-curves $\tilde{E}_j$. Their intersection numbers
do not change under this map, so none of them is a $(-1)$-curve.
The canonical divisor of $Y$ is 
$K_Y=\sum \tilde{L}_i + \sum F_j$. This shows that 
the surface $Y$ is minimal. Its Chern numbers are
$$c_1^2(Y)=k-n+n_0, \quad 
c_2(Y)=48+2(k-n)-n_0.$$

Now let us evaluate the known conditions on
$c_1^2(Y)$ and $c_2(Y)$.

First of all there is Noether's formula \cite{BPV}, p.20,
$$12 \chi(\OO_Y)=c_1^2(Y)+c_2(Y).$$
It shows that
$$c_1^2(Y)+c_2(Y) = 48+3(k-n)$$
is divisible by $12$. This is just the modulo-$4$ condition
$(\lambda 6)$ from section 2, nothing new.

Then there is the famous Miyaoka-Yau-inequality \cite{BPV}, p. 212,
$$c_1^2(Y) \leq 3c_2(Y).$$
In our case it reads
\begin{eqnarray*}
k-n+n_0 & \leq & 144+6(k-n)-3n_0, \\
5n+4n_0 & \leq & 144+5k.
\end{eqnarray*}
In our range ($n \leq 10, n_0 \leq 10$) this is not a surprise.

Finally consider Noether's inequality \cite{BPV}, p.211,
$$5c_1^2(Y)-c_2(Y)+36 \geq 0$$
for minimal surfaces of general type. In our case it reads
$$5(k-n)+5n_0-(48+2(k-n)-n_0)+36
= 3(k-n)+6n_0-12 \geq 0.$$
Leaving aside the arrangements $\Lambda(6), \, \Lambda_1(8),
\, \Lambda_2(8), \Lambda_1(10)$ and $\Lambda_2(10)$,
which lead to elliptic or abelian surfaces $Y$, 
we find that this inequality always holds. Although, for the arrangements
$\Lambda_5(10),...,\Lambda_9(10)$ it is an equality, the Chern numbers being
$$c_1^2(Y)=4, \quad c_2(Y)=56.$$

\VV

\end{document}